\documentclass[12pt]{article}
\usepackage{amssymb,amsfonts,amsmath}\usepackage{graphicx}\usepackage{amsfonts}

\def\0{\leqno}
\def\({\left(}
\def\){\right)}
\def\<{\left<}
\def\>{\right>}
\def\ov{\overline}
\def\4{\subseteq }
\def\dd{\displaystyle}

\def\bit{\begin{itemize}}
\def\eit{\end{itemize}}
\def\barr{\begin{array}}
\def\earr{\end{array}}

\def\X#1#2{\stb{#1}{#2}{\mbox{\Huge$\times$}}}
\def\Z{{\rlap{$\kern2pt{\rm Z}$}{\rm Z}\,}}
\def\bld#1#2{{\buildrel{#1}\over{#2}}}
\def\st#1#2{{\mathrel{\mathop{#2}\limits_{#1}}{}\!}}
\def\stb#1#2#3{{\st{{#1}}{\bld{{#2}}{#3}}{}\!}}

\def\a{{\alpha}}

\def\dd{\displaystyle}
\title{\bf An arithmetic method of counting the subgroups of a finite abelian group}
\author{ Marius T\u arn\u auceanu}
\date{October 1, 2010}

\begin{document}

\maketitle

\begin{abstract}
The main goal of this paper is to apply the arithmetic method
developed in our previous paper \cite{13} to determine the number
of some types of subgroups of finite abelian groups.
\end{abstract}



\section{Introduction}

One of the most important problems of the  combinatorial abelian
group theory is to determine the number of subgroups of a finite
abelian group. This topic has enjoyed a constant evolution
starting with the first half of the $20^{\mathrm{th}}$ century. Since a finite
abelian group is a direct product of abelian $p$-groups, the above
counting problem is reduced to $p$-groups. Formulas which give the
number of subgroups of type $\mu$ of a finite $p$-group of type
$\lambda$ were established by S. Delsarte (see \cite{7}), P.E.
Djubjuk (see \cite{8}) and Y. Yeh (see \cite{15}). An excellent
survey on this subject together with connections to symmetric
functions was written by M.L. Butler (see \cite{5}) in 1994.
Another way to find the total number of subgroups of finite
abelian $p$-groups is presented in \cite{6} and applied for rank
two $p$-groups, as well as for elementary abelian $p$-groups.
Also, remind here the paper \cite{1} which gives an explicit
formula for the number of subgroups in a finite abelian $p$-group
by using divisor functions of matrices.

The starting point for our discussion is given by the paper
\cite{13} (see also Section I.2 of \cite{14}), where we introduced
and studied the concept of fundamental group lattice, that is the
subgroup lattice of a finite abelian group. These lattices were
successfully used to solve the problem of existence and uniqueness
of a finite abelian group whose subgroup lattice is isomorphic to
a fixed lattice (see Proposition 2.8 (\S\ 2.1) of \cite{13}). Some
steps in finding the total number of subgroups of several
particular finite abelian groups have been made in Section 2.2 of
\cite{13}, too.

The purpose of the current paper is to extend the above study, by
applying the fundamental group lattices in counting some different
types of subgroups of finite abelian groups. Explicit formulas are
obtained for the number of subgroups of a given order in a finite
abelian $p$-group of rank 2, improving Proposition 2.9 (\S\ 2.2)
of \cite{13}, and for the number of maximal subgroups and cyclic
subgroups of a given order of {\it arbitrary} finite abelian
groups. The number of elements of a prescribed order in such a
group will be also found.

The paper is organized as follows. In Section 2 we recall the
notion of fundamental group lattice and its basic properties.
Section 3 deals with the number of subgroups of finite abelian
groups. In Section 4 the precise expressions for the number of
cyclic subgroups, as well as for the number of elements of a given
order in a finite abelian group will be determined. In the final
section some conclusions and further research directions are
indicated.

Most of our notation is standard and will usually not be  repeated
here. Basic definitions and results on lattices (respectively on
groups) can be found in \cite{9} (respectively in \cite{12}). For
subgroup lattice concepts we refer the reader to \cite{10} and
\cite{14}.

\section{Fundamental group lattices}

Let $G$ be an abelian group of order $n$ and $L(G)$ be the
subgroup lattice of $G$. By the fundamental theorem of finitely
generated abelian groups, there exist (uniquely determined by $G$)
the numbers $k\in{\rm I\!N}^*$, $d_1,d_2,...,d_k\in{\rm
I\!N}\setminus\{0,1\}$ satisfying $d_1|d_2|...|\hspace{0,5mm}d_k$,
$d_1d_2\cdots d_k =n$ and $$G\cong\X{i=1}{k}\Z_{d_i}.$$This
decomposition of a finite abelian group into a direct product of
cyclic groups together with the form of subgroups of $\Z^k$ (see
Lemma 2.1, \S\ 2.1, \cite{13}) leads us to the concept of
fundamental group lattice, defined in the following manner:

 Let $k{\ge}1$ be a natural number. Then, for each
 $(d_1,d_2,...,d_k)\in({\rm I\!N}{\setminus}\{0,1\})^k,$ we consider the set
 $L_{(k;d_1,d_2,...,d_k)}$ consisting of all matrices
 $A=(a_{ij})\in{\cal M}_k(\Z)$ which have the following properties:\bigskip

 $\begin{tabular}{rl}
i)& $a_{ij}=0,$ for any $i>j,$\vspace*{1,2mm}\\
ii)& $0\le a_{1j},a_{2j},...,a_{j-1j}<a_{jj},$ for any
 $j=\overline{1,k}$,\vspace*{1,2mm}\\
iii)& 1)\ $a_{11}|d_1,$\vspace*{1,2mm}\\
 &2)\ $a_{22}|\Bigl(d_2,d_1\,\displaystyle\frac{a_{12}}{a_{11}}\Bigr),$\end{tabular}$

$\begin{tabular}{rl}
 \ \ \ \ \ \ &3)\
 $a_{33}|\Bigl(d_3,d_2\,\displaystyle\frac{a_{23}}{a_{22}},d_1\,\displaystyle\frac{\left|\begin{array}{ll}a_{12}&a_{13}\\
 a_{22}&a_{23}\end{array}\right|}{a_{22}a_{11}}\Bigr),$\\
 &$\vdots$\\
 &k)\ $a_{kk}|\Bigl(d_k,d_{k-1}\,\displaystyle\frac{a_{k-1k}}{a_{k-1\,k-1}},d_{k-2}\,
 \displaystyle\frac{\left|\begin{array}{ll}a_{k-2\,k-1}&a_{k-2k}\\
 a_{k-1\,k-1}&a_{k-1k}\end{array}\right|}
 {a_{k-1\,k-1}a_{k-2\,k-2}},...,$\vspace*{4,5mm}\\
 &\hfill$d_1\,\displaystyle\frac{\left|\begin{array}{llcl}
 a_{12}&a_{13}&\cdots&a_{1k}\\
 a_{22}&a_{23}&\cdots&a_{2k}\\
 \vdots&\vdots&&\vdots\vspace*{-3mm}\\
 0&0&\cdots&a_{k-1\,k}\end{array}\right|}{a_{k-1\,k-1}a_{k-2\,k-2}...
 a_{11}}\Bigr),$
 \end{tabular}$
 \bigskip

 \noindent
 where by $(x_1,x_2,...,x_m)$ we denote the greatest common divisor
 of the numbers $x_1,x_2,...,x_m\in\Z$.
 On the set $L_{(k;d_1,d_2,...,d_k)}$ we introduce the next
 partial ordering relation (denoted by $\leq $), as follows: for
 $ A=(a_{ij}),$ $ B=(b_{ij})\in L_{(k;d_1,d_2,...,d_k)},$ put
 $A\le B$ if and only if the relations\vspace*{-0,1mm}
 \begin{itemize}
 \item[$1)'$] $b_{11}|a_{11},$\vspace*{-1,5mm}
 \item[$2)'$] $b_{22}|\Bigl(a_{22},\displaystyle\frac{\left|\begin{array}{ll}
 a_{11}&a_{12}\\
 b_{11}&b_{12}
 \end{array}\right|}{b_{11}}\Bigr),$\vspace*{-1,5mm}
 \item[$3)'$] $b_{33}|\Bigl(a_{33},\displaystyle\frac{\left|\begin{array}{ll}
 a_{22}&a_{23}\\
 b_{22}&b_{23}\\
 \end{array}\right|}{b_{22}},
 \displaystyle\frac{\left|\begin{array}{lll}
 a_{11}&a_{12}&a_{13}\\
 b_{11}&b_{12}&b_{13}\\
 0&b_{22}&b_{23}
 \end{array}\right|}{b_{22}b_{11}}\Bigr),$\vspace*{-1,5mm}
 \item[$\vdots$]\vspace*{-3,5mm}
 \item[${\rm k})'$] $b_{kk}|\Bigl(a_{kk},\displaystyle\frac{\left|\begin{array}{ll}
 a_{k-1\,k-1}&a_{k-1\,k}\\
 b_{k-1\,k-1}&b_{k-1\,k}
 \end{array}\right|}{b_{k-1\,k-1}},
 \displaystyle\frac{\left|\begin{array}{lll}
 a_{k-2\,k-2}&a_{k-2\,k-1}&a_{k-2\,k}\\
 b_{k-2\,k-2}&b_{k-2\,k-1}&b_{k-2\,k}\\
 0&b_{k-1\,k-1}&b_{k-1\,k}
 \end{array}\right|}{b_{k-1\,k-1}b_{k-2\,k-2}},...,$\bigskip\\
 \ \hspace*{6cm}$\displaystyle\frac{\left|\begin{array}{llcl}
 a_{11}&a_{12}&\cdots&a_{1k}\\
 b_{11}&b_{12}&\cdots&b_{1k}\\
 \vdots&\vdots&&\vdots\\
 0&0&\cdots&b_{k-1\,k}\end{array}\right|}{b_{k-1\,k-1}b_{k-2\,k-2}...b_{11}}\Bigr)$\vspace*{-1,5mm}
 \end{itemize}
\bigskip

 \noindent hold.  Then $(L_{(k;d_1,d_2,...,d_k)},\le)$ is a complete  modular
 lattice, which is called    a {\it   fundamental group lattice of degree
 $k$}.

 Moreover, from Proposition 2.2, \S\ 2.1, \cite{13}, we
 know that
 $L_{(k;d_1,d_2,...,d_k)}$
 is isomorphic to $L(G)$ and so the problem of counting the subgroups of
 $G$  can be translated into an arithmetic problem: finding the number of elements of $L_{(k;d_1,d_2,...,d_k)}$.

 On the other hand, if $n=p_1^{n_1}p_2^{n_2}...p_m^{n_m}$ is the
 decomposition of $n$ as a product of prime factors and
 $$G\cong\X{i=1}mG_i$$is the corresponding primary decomposition of $G$,
 then it is well-known that we have
 $$L(G)\cong\X{i=1}mL(G_i).$$
 The above lattice isomorphism shows that
 $$|L(G)|=\displaystyle\prod_{i=1}^m|L(G_i)|,$$therefore our counting problem is reduced to $p$-groups. In this way, we
need to investigate only fundamental group lattices of type
$L_{(k;p^{\a_1},p^{\a_2},...,p^{\a_k})}$, where $p$ is a prime and
$1\le\a_1\le\a_2\le...\le\a_k$.\newpage
 \noindent Concerning these lattices, the following elementary remarks will be very useful:\vspace*{-1mm}
\begin{itemize}
\item[a)] The order of the subgroup of $\X{i=1}k\Z_{p^{\a_i}}$ corresponding to the matrix\medskip\break $A=(a_{ij})\in L_{(k;p^{\a_1},p^{\a_2},...,p^{\a_k})}$ is $$\dd\frac{p^{\hspace{1mm}\sum\limits_{i=1}^k\a_i}}{\prod\limits_{i=1}^ka_{ii}}\,\cdot$$
\item[b)] The subgroup of $\X{i=1}k\Z_{p^{\a_i}}$ corresponding to the matrix\break $A=(a_{ij})\in L_{(k;p^{\a_1},p^{\a_2},...,p^{\a_k})}$  is cyclic if and only if
 $<(\overline 0^1,\overline0^2,...,\bar a^k_{kk})>\ \subseteq\smallskip\break
<(\overline 0^1,\overline0^2,...,\bar a^{k-1}_{k-1\,k-1},\bar a^k_{k-1\,k})>\ \subseteq\ \cdots\ \subseteq\
<(\bar a^1_{11}\bar a^2_{12},...,\bar a^k_{1k})>,$  where, for\smallskip\break every $i=\overline{1,k}$, we denote by $\bar x^i$ the image of an element $x\in\Z$ through the canonical homomorphism: $\Z\to\Z_{p^{\a_i}}.$
\item[c)] If $A=(a_{ij})$ is an element of $L_{(k;p^{\a_1},p^{\a_2},...,p^{\a_k})}$, then the linear system
$$\barr{l}
A^\top\\\ \\\ \\\ \earr\!\!\!\!\left(\barr{c}x_1\\ x_2\\\vdots\\ x_k\earr\right)=\left(\barr{c}p^{\a_1}\\ p^{\a_2}\\\vdots\\ p^{\a_k}\earr\right)$$admits solutions in $\Z^k$.\vspace*{-1mm}\end{itemize}

\section{The number of subgroups of a finite abelian group}

As we have seen in the previous section, in order to determine the
number of  subgroups of finite abelian groups it suffices to
reduce the study to $p$-groups and our problem is equivalent to
the counting of elements of the fundamental group lattice
$L_{(k;p^{\a_1},p^{\a_2},...,p^{\a_k})}$. This consists of all
matrices of integers $A=(a_{ij})_{i,j=\overline{1,k}}$ satisfying
the conditions:

$$\left\{\begin{tabular}{rl}
i)& $a_{ij}=0,$ for any $i>j,$\vspace*{3mm}\\
ii)& $0\le a_{1j},a_{2j},...,a_{j-1j}<a_{jj},$ for any
 $j=\overline{1,k}$,\vspace*{3mm}\\
iii)& 1)\ $a_{11}|p^{\a_1},$\vspace*{2mm}\\
 &2)\ $a_{22}|\Bigl(p^{\a_2},p^{\a_1}\,\displaystyle\frac{a_{12}}{a_{11}}\Bigr),$\vspace*{2mm}\\
 &3)\
 $a_{33}|\Bigl(p^{\a_3},p^{\a_2}\,\displaystyle\frac{a_{23}}{a_{22}},p^{\a_1}\,\displaystyle\frac{\left|\begin{array}{ll}a_{12}&a_{13}\\
 a_{22}&a_{23}\end{array}\right|}{a_{22}a_{11}}\Bigr),$\\
 &$\vdots$\\
 &k)\ $a_{kk}|\Bigl(p^{\a_k},p^{\a_{k-1}}\,\displaystyle\frac{a_{k-1k}}{a_{k-1\,k-1}},p^{\a_{k-2}}\,
 \displaystyle\frac{\left|\begin{array}{ll}a_{k-2\,k-1}&a_{k-2k}\\
 a_{k-1\,k-1}&a_{k-1k}\end{array}\right|}
 {a_{k-1\,k-1}a_{k-2\,k-2}},...,$\vspace*{4,9mm}\\
 &\hfill$p^{\a_1}\,\displaystyle\frac{\left|\begin{array}{llcl}
 a_{12}&a_{13}&\cdots&a_{1k}\\
 a_{22}&a_{23}&\cdots&a_{2k}\\
 \vdots&\vdots&&\vdots\\
 0&0&\cdots&a_{k-1\,k}\end{array}\right|}{a_{k-1\,k-1}a_{k-2\,k-2}...
 a_{11}}\Bigr).$
 \end{tabular}\right.\leqno(*)$$
 \medskip

\noindent An explicit formula for
$|L_{(k;p^{\a_1},p^{\a_2},...,p^{\a_k})}|$, and  consequently for
$|L(\X{i=1}k\Z_{p^{\a_i}})|$, can be easily obtained in the
particular case $\a_1=\a_2=\cdots=\a_k=1$ (see Proposition 2.12,
\S\ 2.2, \cite{13}).

\bigskip\noindent{\bf Proposition 3.1.} {\it For
$\a\in\{0,1,...,k\}$, the number of all subgroups of order
$p^{k-\a}$ in the finite elementary abelian $p$-group $\Z^k_p$ is
$1$ if $\a=0$ or $\a=k$, and $\displaystyle\sum_{1\le
i_1<i_2<...<i_\a\le k}p^{i_1+i_2+...+i_\a-\frac{\a(\a+1)}2}$ if
$1\le\a\le k-1.$ In particular, the total number of subgroups of
$\Z^k_p$ is $2+\displaystyle\sum_{\a=1}^{k-1}
\displaystyle\sum_{1\le i_1<i_2<...<i_\a\le
k}p^{i_1+i_2+...+i_\a-\frac{\a(\a+1)}2}$.}\bigskip\medskip

In the general case, our method gives an immediate result in
counting the maximal subgroups of $\X{i=1}k\Z_{p^{\a_i}}.$ By the
first remark of Section 2, such a subgroup corresponds to a matrix
$A=(a_{ij})\in L_{(k;p^{\a_1},p^{\a_2},...,p^{\a_k})}$ satisfying
$\displaystyle\prod_{i=1}^ka_{ii}=p$. Then $a_{ii}=p$ for some
$i\in\{1,2,...,k\}$ and $a_{jj}=1$ for all $j\ne i$. From the
condition ii) of $(*)$ we get $a_{1j}=a_{2j}=\cdots=a_{j-1j}=0$,
for all $j\ne i$. Remark also that the condition iii) is satisfied
and thus the elements $a_{1i},a_{2i},...,a_{i-1i}$ can be chosen
arbitrarily from the set $\{0,1,...,p-1\}$. Therefore we have
$p^{i-1}$ distinct solutions of the system $(*)$. Summing up these
quantities for $i=\overline{1,k}$, we determine the number of
maximal subgroups of the finite abelian $p$-group
$\X{i=1}k\Z_{p^{\a_i}}.$

\bigskip\noindent{\bf Proposition 3.2.} {\it The number of maximal
subgroups of $\X{i=1}k\Z_{p^{\a_i}}$ is
$\displaystyle\frac{p^k-1}{p-1}\,\cdot$}\bigskip\medskip

Next, we return to the problem of finding the total number  of
subgroups of $\X{i=1}k\Z_{p^{\a_i}}$. We shall apply our method
for rank two abelian $p$-groups, i.e. when $k=2$ (clearly, it can
be extended in a natural way for an arbitrary $k$). Note also that
the following theorem improves Proposition 2.9, \S\ 2.2,
\cite{13}, by indicating the number of subgroups of a fixed order
in such a group and by giving a proof founded on fundamental group
lattices.

\bigskip\noindent{\bf Theorem 3.3.} {\it  For every $0\le\a\le\a_1+\a_2$,  the number of all subgroups of order $p^{\a_1+\a_2-\a}$ in the finite abelian $p$-group $\Z_{p^{\a_1}}\times\Z_{p^{\a_2}}$ is:

$$\left\{\barr{lll}
\displaystyle\frac{p^{\a+1}-1}{p-1}\,,&\mbox{ if }&0\le\a\le\a_1\vspace*{3mm}\\
\displaystyle\frac{p^{\a_1+1}-1}{p-1}\,,&\mbox{ if }&\a_1\le\a\le\a_2\vspace*{3mm}\\
\displaystyle\frac{p^{\a_1+\a_2-\a+1}-1}{p-1}\,,&\mbox{ if
}&\a_2\le\a\le\a_1+\a_2.\earr\right.$$In particular, the total
number of subgroups of $\Z_{p^{\a_1}}\times\Z_{p^{\a_2}}$ is
$$\frac1{(p{-}1)^2}\left[(\a_2{-}\a_1{+}1)p^{\a_1{+}2}{-}(\a_2{-}\a_1{-}1)p^{\a_1{+}1}{-}(\a_1{+}\a_2{+}3)p{+}(\a_1{+}\a_2+1)\right].$$}\bigskip

\noindent{\bf Proof.} Let $A=(a_{ij})$ be a solution of $(*)$ for
$k=2$,  corresponding to a subgroup of order $p^{\a_1+\a_2-\a}$.
In this situation, the condition iii) of $(*)$ becomes

$$a_{11}|p^{\a_1}\mbox{\ \ and\ \ }a_{22}|\left(p^{\a_2},p^{\a_1}\,\frac{a_{12}}{a_{11}}\right).$$

\noindent Put $a_{11}=p^i,$ where $0\le i\le\a_1.$ Then $a_{22}=p^{\a-i}$ and so $p^{\a-i}|(p^{\a_2},p^{\a_1-i}a_{12})$, that is $p^{\a-i}|p^{\a_1-i}(p^{\a_2-\a_1+i},a_{12})$. If $0\le\a\le\a_1$, we must have $i\le\a$ and the above condition is satisfied by all $a_{12}<p^{\a-i}.$ So, one obtains $p^{\a-i}$ distinct solutions of $(*)$, which implies that the number of subgroups of order $p^{\a_1+\a_2-\a}$ in $\Z_{p^{\a_1}}\times\Z_{p^{\a_2}}$ is in this case
$$S_1(\a)=\sum_{i=0}^\a p^{\a-i}=\frac{p^{\a+1}-1}{p-1}\,\cdot\0(1)$$
Suppose now that $\a_1\le\a\le\a_2$. Then
$p^{\a_1-\a}|(p^{\a_2-\a_1+i},a_{12})$  and thus $a_{12}$ can be
any multiple of $p^{\a_1-\a}$ in the set $\{0,1,...,p^{\a-i}-1\}$.
It results $p^{\a_1-i}$ distinct solutions of $(*)$ and the number
of subgroups of order $p^{\a_1+\a_2-\a}$ in
$\Z_{p^{\a_1}}\times\Z_{p^{\a_2}}$ is in this case
$$S_2(\a)=\sum_{i=0}^{\a_1}p^{\a_1-i}=\frac{p^{\a_1+1}-1}{p-1}\,\cdot\0(2)$$
Finally, assume that $\a_2\le\a\le\a_1+\a_2$. We must have
$\a_1-\a\le\a_2-\a_1+i$  and the number of distinct solutions of
$(*)$ is again $p^{\a_1-i}.$ Thus the number of subgroups of order
$p^{\a_1+\a_2-\a}$ in $\Z_{p^{\a_1}}\times\Z_{p^{\a_2}}$ is in
this case
$$S_2(\a)=\sum_{i=\a-\a_2}^{\a_1}p^{\a_1-i}=\frac{p^{\a_1+\a_2-\a+1}-1}{p-1}\,\cdot\0(3)$$
By using the equalities (1), (2) and (3), one obtains the total
number of  subgroups of $\Z_{p^{\a_1}}\times\Z_{p^{\a_2}}$, namely
$$\barr{c}
\displaystyle\sum_{\a=0}^{\a_1}S_1(\a)+\sum_{\a=\a_1+1}^{\a_2}S_2(\a)+\sum_{\a=\a_2+1}^{\a_1+\a_2}S_3(\a)=\displaystyle\frac1{(p-1)^2}\left[(\a_2-\a_1+1)p^{\a_1+2}-\right.\vspace*{2,5mm}\\
-\left.(\a_2-\a_1-1)p^{\a_1+1}-(\a_1+\a_2+3)p+(\a_1+\a_2+1)\right],
\earr$$
which completes our proof. $\scriptstyle\Box$\bigskip

In the following let us denote by $f_p(i,j)$ the number of all
subgroups of the  finite abelian $p$-group
$\Z_{p^i}\times\Z_{p^j}$ $(i\le j)$, determined in Theorem 3.3.
Note that we have
$$\barr{lcl}
f_p(i,j)=\vspace*{1,5mm}\\
=\displaystyle\frac1{(p-1)^2}\left[(j-i+1)p^{i+2}-(j-i-1)p^{i+1}-(i+j+3)p+(i+j+1)\right]=\vspace*{1,5mm}\\
=(j-i+1)p^i+(j-i+3)p^{i-1}+\cdots+(i+j-1)p+(i+j+1). \earr$$ Put
$f_p(i,j)=f_p(j,i)$, for all $i>j$, and let $n$ be a fixed
positive integer and $A_p(n)$ be the matrix
$(f_p(i,j))_{i,j=\overline{0,n}}.$ Then $A_p(n)$ induces a
quadratic form $\sum\limits_{i,j=0}^nf_p(i,j)X^iY^j$. Because
$$\det A_p(n){=}(p{-}1)p^{n-1}\det A_p(n{-}1),$$by induction on $n$
one easily obtains $$\det A_p(n)=(p-1)^np^{\frac{n(n-1)}2}, \mbox{
for any } n\ge1.$$Hence, we have proved the next two corollaries.

\bigskip\noindent{\bf Corollary 3.4.} {\it The quadratic form
$\sum\limits_{i,j=0}^nf_p(i,j)X^iY^j$ induced by the matrix $A_p(n)$ is positive definite, for all $n\in{\rm I\!N}^*.$}

\bigskip\noindent{\bf Corollary 3.5.} {\it All eigenvalues of the
matrix $A_p(n)$ are positive, for all $n\in{\rm I\!N}^*$.}

\section{The number of cyclic subgroups of a finite abelian group}

Another interesting application of fundamental group lattices (not
studied in \cite{13}) is the counting of  cyclic subgroups of
finite abelian groups. First of all, we obtain this number for a
finite abelian $p$-group of rank 2. By the second remark of
Section 2, the subgroup of $\Z_{p^{\a_1}}\times\Z_{p^{\a_2}}$
determined by the matrix $A=(a_{ij})$ is cyclic if and only if
$<(\overline0^1,\bar a^2_{22})>\ \subseteq\ <(\bar a^1_{11},\bar
a^2_{12})>.$ This necessary and sufficient condition can be
rewritten in the following manner.

\bigskip\noindent{\bf Lemma 4.1.} {\it The subgroup of
$\Z_{p^{\a_1}}\times\Z_{p^{\a_2}}$  corresponding to the matrix
$A=(a_{ij})\in L_{(2;p^{\a_1},p^{\a_2})}$ is cyclic if and only if
$a_{22}=\left(p^{\a_2},p^{\a_1}\,\displaystyle\frac{a_{12}}{a_{11}}\right).$}\bigskip

\noindent{\bf Proof.} If $<(\ov0^1,\bar a^2_{22})>\ \4\ <(\bar
a^1_{11},\bar a^2_{12})>,$
 then we can choose an integer $x$ such that $(\ov0^1,\bar a^2_{22})=x(\bar a^1_{11},\bar
 a^2_{12}).$
 It results $p^{\a_1}|xa_{11}$ and $p^{\a_2}|xa_{12}-a_{22},$ therefore there
 exist
 $y,z\in\Z$ satisfying $xa_{11}=yp^{\a_1}$ and
 $xa_{12}-a_{22}=zp^{\a_2}$.
 These equalities imply that $a_{22}=-zp^{\a_2}+yp^{\a_1}\,\displaystyle
 \frac{a_{12}}{a_{11}},$
  which together with the condition 2) of iii)
  in $(*)$ show that $a_{22}=\(p^{\a_2},p^{\a_1}\,\displaystyle\frac{a_{12}}{a_{11}}\)\cdot$

Conversely, suppose that
$a_{22}=\(p^{\a_2},p^{\a_1}\,\displaystyle\frac{a_{12}}{a_{11}}\)\cdot$
Then there are $y,z\in \Z$ with
$a_{22}=-zp^{\a_2}+yp^{\a_1}\,\displaystyle
\frac{a_{12}}{a_{11}}\,\cdot$ Taking
$x=y\,\dd\frac{p^{\a_1}}{a_{11}}\in\Z$, we  easily obtain
$(\ov0^1,\bar a^2_{22})=x(\bar a^1_{11},\bar a^2_{12})$ and so
$<(\ov0^1,\bar a^2_{22})>$ is contained in $<(\bar a^1_{11},\bar
a^2_{12})>.$ $\scriptstyle\Box$\bigskip

By using the above lemma, the problem of finding the number of
cyclic  subgroups of $\Z_{p^{\a_1}}\times\Z_{p^{\a_2}}$ reduces to
an elementary arithmetic exercise.

\bigskip\noindent{\bf Theorem 4.2.} {\it For every
$0\le\a\le\a_2$, the number of cyclic subgroups of order $p^\a$ in
the finite abelian $p$-group $\Z_{p^{\a_1}}\times\Z_{p^{\a_2}}$
is:
$$\left\{\barr{lll}
1,&\mbox{ if }&\a=0\vspace*{1,5mm}\\
p^\a+p^{\a-1},&\mbox{ if }&1\le\a\le\a_1\vspace*{1,5mm}\\
p^{\a_1},&\mbox{ if }&\a_1<\a\le\a_2.\earr\right.$$
In particular, the number of all cyclic subgroups of $\Z_{p^{\a_1}}\times\Z_{p^{\a_2}}$ is
$$2+2p+\cdots+2p^{\a_1-1}+(\a_2-\a_1+1)p^{\a_1}.$$}

\noindent{\bf Proof.} Denote by $g^2_p(\a)$ the number of  cyclic
subgroups of order $p^\a$ in
\mbox{$\Z_{p^{\a_1}}{\times}\Z_{p^{\a_2}}$} and let
$A{=}(a_{ij})\in L_{(2;p^{\a_1},p^{\a_2})}$ be the matrix
corresponding to such a subgroup. Then $a_{11}|p^{\a_1},$
$a_{22}{=}\(p^{\a_2},p^{\a_1}\,\dd\frac{a_{12}}{a_{11}}\)$ and
$a_{11}a_{22}{=}p^{\a_1+\a_2-\a}.$ Taking $a_{11}{=}p^i$ with
$0\le i\le\a_1$, we obtain
$$a_{22}=p^{\a_1+\a_2-\a-i}=(p^{\a_2},p^{\a_1-i}a_{12})=p^{\a_1-i}(p^{\a_2-\a_1+i},a_{12}),$$which implies that
$$p^{\a_2-\a}=\(p^{\a_2-\a_1+i},a_{12}\).\0(4)$$
Clearly, for $\a=0$ it results $a_{11}=p^{\a_1},$
$a_{22}=p^{\a_2},$ $a_{12}=0$, and thus
$$g^2_p(0)=1.\0(5)$$
For $1\le\a\le\a_1$ we must have $\a_1-\a\le i.$ If $i=\a_1-\a,$
the condition  (4) is equivalent to $p^{\a_2-\a}|a_{12}$,
therefore $a_{12}$ can be chosen in $p^\a$ ways. If $\a_1-\a+1\le
i$, (4) is equivalent to
$$p^{\a_2-\a}|a_{12}\mbox{\ \ and\ \ }p^{\a_2-\a+1}\nmid a_{12}.$$
There are $p^{\a_1-i}\hspace{-0,5mm}-p^{\a_1-i-1}$ elements of the
set $\{0,1,...,p^{\a_1+\a_2-\a-i}\}$ which satisfy the previous
relations. So, one obtains
$$g^2_p(\a)=p^\a+\sum_{i=\a_1-\a+1}^{\a_1}\(p^{\a_1-i}-p^{\a_1-i-1}\)=p^\a+p^{\a-1},\mbox{ for }1\le\a\le\a_1.\0(6)$$
Mention that if $\a_1<\a\le\a_2$, then the condition $\a_1-\a\le
i$ is satisfied by all $i=\ov{1,\a_1}$, and hence
$$g^2_p(\a)=\sum_{i=0}^{\a_1}\(p^{\a_1-i}-p^{\a_1-i-1}\)=p^{\a_1},\mbox{ for }\a_1<\a\le\a_2.\0(7)$$
Now, the equalities (5)-(7) give us the total number of cyclic
subgroups of $\Z_{p^{\a_1}}\times\Z_{p^{\a_2}}$, namely
$$\hspace{-50mm}1+\dd\sum_{\a=1}^{\a_1}\(p^\a+p^{\a-1}\)+\dd\sum_{\a=\a_1+1}^{\a_2}p^{\a_1}=$$
$$\hspace{22mm}=\dd\frac1{p-1}
\left[(\a_2-\a_1+1)p^{\a_1+1}-(\a_2-\a_1-1)p^{\a_1}-2\right]=$$
$$=2+2p+\cdots+2p^{\a_1-1}+(\a_2-\a_1+1)p^{\a_1}$$and
our proof is  finished. $\scriptstyle\Box$
\bigskip

The above method can be used for an arbitrary $k>2$, too. In order to do this we need to remark that
$$g^2_p(\a)=\frac{p^\a h^1_p(\a)-p^{\a-1}h^1_p(\a-1)}{p^\a-p^{\a-1}}\,,\mbox{ for all }\a\ne0,$$where
$$h^1_p(\a)=\left\{\barr{lll}
p^\a,&\mbox{ if }&0\le\a\le\a_1\vspace*{1,5mm}\\
p^{\a_1},&\mbox{ if }&\a_1\le\a.\earr\right.$$\medskip

\noindent This equality extends to the general case in the
following way.

\bigskip\noindent{\bf Theorem  4.3.} {\it For every $1\le\a\le\a_k,$ the
number of cyclic subgroups of order $p^\a$ in the finite abelian
$p$-group $\X{i=1}k\Z_{p^{\a_i}}$ is
$$g^k_p(\a)=\frac{p^\a h^{k-1}_p(\a)-p^{\a-1}h^{k-1}_p(\a-1)}{p^\a-p^{\a-1}}\,,$$
where
$$h^{k-1}_p(\a)=\left\{\barr{lll}
p^{(k-1)\a},&\mbox{ if }&0\le\a\le\a_1\vspace*{1,5mm}\\
p^{(k-2)\a+\a_1},&\mbox{ if }&\a_1\le\a\le\a_2\\
\vdots\\
p^{\a_1+\a_2+...+\a_{k-1}},&\mbox{ if }&\a_{k-1}\le\a.\earr\right.$$}

Note that $g^k_p(0)=1$ and the number of all cyclic subgroups of
$\X{i=1}k\Z_{p^{\a_i}}$  can be easily determined from Theorem
4.3. Since the numbers of cyclic subgroups and of elements of a
given order in a finite abelian $p$-group are closely connected
(through the well-known Euler's function $\varphi$), we also infer
the following consequence of Theorem 4.3.

\bigskip\noindent{\bf Corollary  4.4.} {\it The number of all elements of order $p^\a$, $1\le\a\le\a_k$,
 in the finite abelian $p$-group $\X{i=1}k\Z_{p^{\a_i}}$ is
$$g^k_p(\a)\varphi(p^\a)=g^k_p(\a)\(p^\a-p^{\a-1}\)=p^\a h^{k-1}_p(\a)-p^{\a-1}h^{k-1}_p(\a-1).$$}

As we have seen in Section 2, counting the subgroups of finite
abelian groups can be reduced to $p$-groups. The same thing can be
also said for cyclic subgroups and for elements of a given order
in an arbitrary finite abelian group $G$. Suppose that
$p_1^{n_1}p_2^{n_2}...p_m^{n_m}$ is the decomposition of $|G|$ as
a product of prime factors and let $\X{i=1}mG_i$ be the
corresponding primary decomposition of $G$. Then every cyclic
subgroup $H$ of order $p_1^{\a_1}p_2^{\a_2}...p_m^{a_m}$ of $G$
can be uniquely written as a direct product $\X{i=1}mH_i$, where
$H_i$ is a cyclic subgroup of order $p_i^{\a_i}$ of $G_i$,
$i=\ov{1,m}$. This remark leads to the following result, that
generalizes Theorem 4.3. and Corollary 4.4.

\bigskip\noindent{\bf Corollary 4.5.} {\it Under the previous
hypotheses, for every $(\a_1,\a_2,...,\a_m)\in{\rm I\!N}^m$ with
$\a_i\le n_i,$ $i=\ov{1,m},$ the number of cyclic subgroups
$($respectively of elements$)$ of order $p_1^{\a_1} p_2^{\a_2} ...
p_m^{\a_m}$ in $G$ is
 $$\prod_{i=1}^mg^{k_i}_{p_i}(\a_i)$$
 $($respectively
 $$\prod_{i=1}^mg^{k_i}_{p_i}(\a_i)\varphi(p_i^{\a_i})),$$
 where $k_i$ denotes the number of direct factors of $G_i,\ i=\overline{1,m}.$}

\section{Conclusions and further research}

All our previous results show that the arithmetic method
introduced in \cite{13} and applied in this paper can constitute
an alternative way to study the subgroups of finite abelian
groups. Clearly, it can successfully be used in solving many
computational problems in (finite) abelian group theory. These
will surely be the subject of some further research.
\bigskip

Finally, we mention several open problems concerning this
topic.
\bigskip

\noindent{\bf Problem 5.1.} Extend Theorem 3.3, by indicating
explicit formulas for the number of subgroups of a fixed order and
for the total number of subgroups of a finite abelian $p$-group of
an {\it arbitrary} rank.
\bigskip

\noindent{\bf Problem 5.2.} Let $G$ be a finite abelian group. Use
the description of $L(G)$ given by the above arithmetic method and
the well-known description of $Aut(G)$ to determine the
characteristic subgroups of $G$.
\bigskip

\noindent{\bf Problem 5.3.} By using the defining relations of a
fundamental group lattice, create a computer algebra program that
generates the subgroups of a finite abelian group.
\bigskip

\noindent{\bf Problem 5.4.} Let $G_1$ be a finite abelian group
and $G_2$ be a finite group such that $\mid G_1\mid=\mid
G_2\mid=n$. Denote $\pi_e(G_i)=\{o(a)\mid a\in G_i\}$ and, for
every divisor $d$ of $n$, let $n_i(d)$ be the number of elements
of order $d$ in $G_i$, $i=1,2$ (remark that the numbers $n_1(d)$
are known, by Corollary 4.5). Is it true that the conditions
\begin{itemize}
\item[a)]$\pi_e(G_1)=\pi_e(G_2)$,
\item[b)]$n_1(d)=n_2(d)$, for all $d$,
\end{itemize}
imply the group isomorphism $G_1\cong G_2$? (In other words, study
whether a finite abelian group is determined by the set of its
element orders and by the numbers of elements of any fixed order).
\bigskip

\bigskip

{\bf Acknowledgement.} The authors are grateful to the reviewers
for their remarks which improve the previous version of the paper.

\vspace*{5ex}\small

\hfill
\begin{minipage}[t]{5cm}
Marius T\u arn\u auceanu \\
Faculty of  Mathematics \\
``Al.I. Cuza'' University \\
Ia\c si, Romania \\
e-mail: {\tt tarnauc@uaic.ro}
\end{minipage}


\begin{thebibliography}{10}
\bibitem{1} Bhowmik, G., {\it Evaluation of the divisor function of matrices}, Acta Arith. {\bf 74} (1996), 155-159.
\bibitem{2} Bhowmik, G., Ramar\'e, O., {\it Average orders of multiplicative arithmetical functions of integer matrices}, Acta Arith. {\bf 66} (1994), 45-62.
\bibitem{3} Bhowmik, G., Ramar\'e, O., {\it Algebra of matrix arithmetic}, J. Algebra {\bf 210} (1998), 194-215.
\bibitem{4} Birkhoff, G., {\it Subgroups of abelian groups}, Proc. Lond. Math. Soc. {\bf 38} (1934, 1935), 385-401.
\bibitem{5} Butler, M.L., {\it Subgroup lattices and symmetric functions}, Mem. Amer. Math. Soc., vol. 112, no. {\bf 539}, 1994.
\bibitem{6} C\u alug\u areanu, Gr.G., {\it The total number of subgroups of a finite abelian group}, Sci. Math. Jpn. (1) {\bf 60} (2004), 157-167.
\bibitem{7} Delsarte, S., {\it Functions de M\"obius sur les groupes abeliens finis}, Ann. of Math. {\bf 49} (1948), 600-609.
\bibitem{8} Djubjuk, P.E., {\it On the number of subgroups of a finite abelian group}, Izv. Akad. Nauk SSR Ser. Mat. {\bf 12} (1948), 351-378.
\bibitem{9} Gr\"atzer, G., {\it General lattice theory}, Academic Press, New York, 1978.
\bibitem{10} Schmidt, R., {\it Subgroup lattices of groups}, de Gruyter Expositions in Mathematics 14, de Gruyter, Berlin, 1994.
\bibitem{11} Suzuki, M., {\it On the lattice of subgroups of finite groups}, Trans. Amer. Math. Soc. {\bf 70} (1951), 345-371.
\bibitem{12} Suzuki, M., {\it Group theory}, I, II, Springer Verlag, Berlin, 1982, 1986.
\bibitem{13} T\u arn\u auceanu, M., {\it A new method of proving some classical theorems of abelian groups}, Southeast Asian Bull. Math. (6) {\bf 31} (2007), 1191-1203.
\bibitem{14} T\u arn\u auceanu, M., {\it Groups determined by posets of subgroups}, Ed. Matrix Rom, Bucure\c sti, 2006.
\bibitem{15} Yeh, Y., {\it On prime power abelian groups}, Bull. Amer. Math. Soc. {\bf 54} (1948), 323-327.
\end{thebibliography}
\end{document}